\documentclass[12pt]{amsart}
\advance\hoffset by -0.5 truecm
\usepackage{amssymb}
\newtheorem{Theorem}{Theorem}[section]
\newtheorem{Lemma}[Theorem]{Lemma}

\newtheorem{Remark}[Theorem]{Remark}

\def\V{\mbox{Var}}

\def\Z{{\mathbb Z}}
\def\R\re
\def\V{\bf V}

\def \la{\lambda}

\def \re{{\mathbb R}}

\def \T{{\mathbb T}}

\def \0{\lambda_{0}}
\def \la{\lambda}
\def \ga{\gamma}

\begin{document}
\title[]{Longitudinal KAM-cocycles and action spectra of magnetic flows}

\author[N.S. Dairbekov]{Nurlan S. Dairbekov}
\address{Sobolev Insittute of Mathematics,
Novosibirsk, 630090, Russia}
\email{dair@math.nsc.ru}

\author[G.P. Paternain]{Gabriel P. Paternain}
 \address{ Department of Pure Mathematics and Mathematical Statistics,
University of Cambridge,
Cambridge CB3 0WB, England}
 \email {g.p.paternain@dpmms.cam.ac.uk}




\begin{abstract} Let $M$ be a closed oriented surface and let
$\Omega$ be a non-exact 2-form. Suppose that the magnetic flow $\phi$ of the pair
$(g,\Omega)$ is Anosov. We show that the longitudinal KAM-cocycle of 
$\phi$ is a coboundary if and only the Gaussian curvature is constant and $\Omega$ 
is a constant multiple of the area form thus extending the results in \cite{P2}.
We also show infinitesimal rigidity of the action spectrum of $\phi$ with respect
to variations of $\Omega$. Both results are obtained by showing that
if $G:M\to\mathbb R$ is any smooth function and $\omega$ is any smooth $1$-form on $M$ such that
$G(x)+\omega_{x}(v)$ integrates to zero along any closed orbit of $\phi$, then
$G$ must be identically zero and $\omega$ must be exact.

\end{abstract}

\maketitle

\section{Introduction}

Let $M$ be a closed oriented surface endowed with a Riemannian metric $g$
and let $\Omega$ be a 2-form. The {\it magnetic flow} of the pair 
$(g,\Omega)$ is the flow $\phi$ on the unit sphere bundle $SM$
determined by the equation
\begin{equation}
\frac{D\dot{\gamma}}{dt}=\la(\gamma)\,i\dot{\ga},
\label{maggeo}
\end{equation}
where $i$ indicates rotation by $\pi/2$ according to the orientation
of the surface and $\la$ is the smooth function on $M$ uniquely determined
by $\Omega=\la\Omega_{a}$, where $\Omega_{a}$ is the area
form of $M$. When $\Omega$ vanishes we recover the usual geodesic flow
of the surface. A~curve $\ga:\re\to M$ that solves (\ref{maggeo})
will be called a {\it magnetic geodesic}.

In the present paper we shall study rigidity properties of {\it Anosov}
magnetic flows. The Anosov property means that $T(SM)$
splits as $T(SM)=E^{0}\oplus E^{u}\oplus E^{s}$ in such a way that
there are constants $C>0$ and $0<\rho<1<\eta$ such that $E^{0}$ is
spanned by  the generating vector field of the flow, 
and for all $t>0$ we have
\[\|d\phi_{-t}|_{E^{u}}\|\leq C\,\eta^{-t}\;\;\;\;\mbox{\rm
and}\;\;\;\|d\phi_{t}|_{E^{s}}\|\leq C\,\rho^{t}.\]
The subbundles are then invariant and H\"older continuous and have
smooth integral manifolds, the stable and unstable manifolds,
which define a continuous foliation with smooth leaves.

To any $C^{k}$ volume preserving Anosov flow $\varphi$ on a closed
3-manifold $N$, P. Foulon and B. Hasselblatt \cite{FH} associated
its {\it longitudinal KAM-cocycle}. This is a cocycle that
measures the regularity of the subbundle $E^{u}\oplus E^{s}$ 
The main theorem in \cite{FH} asserts that $E^{u}\oplus E^{s}$ is 
always Zygmund-regular and that the following are
equivalent:
\begin{enumerate}
\item $E^{u}\oplus E^{s}$ is ``little Zygmund'';
\item the longitudinal KAM-cocycle is a coboundary;
\item $E^{u}\oplus E^{s}$ is Lipschitz;
\item $E^{u}\oplus E^{s}$ is $C^{k-1}$;
\item $\varphi$ is a suspension or contact flow.
\end{enumerate}

(A continuous function $f:U\to\re$ on an open set $U\subset \re$ is said to be
{\it Zygmund-regular} if
 $|f(x+h)+f(x-h)-2f(x)|=O(h)$ for all $x$ in $U$. 
The function is said to be {\it ``little Zygmund''} if $|f(x+h)+f(x-h)-2f(x)|=o(h)$.)

It is well known that for flows, a ``choice of time" or equivalently, a choice of speed at which
orbits travel gets reflected on the regularity of the corresponding strong stable and strong unstable distributions.
 The situation is different if we look at the weak unstable and stable bundles $E^{0}\oplus E^{u}$ and
$E^{0}\oplus E^{s}$. S. Hurder and A. Katok proved \cite{HK} that
the weak bundles are always differentiable with Zygmund-regular derivative
and there is a cocycle obstruction to higher regularity given by the first
nonlinear term in the Moser normal form (this explains why Foulon and 
Hasselblatt used the terminology ``longitudinal KAM-cocycle'').

In \cite{P2}, the second author showed that if $\Omega$ is non-exact, $g$ has negative Gaussian curvature $K$
and $\la$ is small enough in the $C^0$ norm, then the longitudinal KAM-cocycle of $\phi$ 
is a coboundary if and only if $K$ and $\la$ are constant.
In the present paper we would like to extend this result to all Anosov magnetic flows, without restrictions
on curvature or the size of $\la$. As shown in \cite{BP} the set of Anosov magnetic flows
can certainly go well beyond small perturbations of Anosov geodesic flows.

\medskip

\noindent {\bf Theorem A.} {\it Let $M$ be a closed oriented
surface endowed with a Riemannian metric $g$ and let $\Omega$ be an arbitrary $2$-form.
Suppose that the magnetic flow $\phi$ of the pair $(g,\Omega)$ is Anosov.
We have:
\begin{enumerate}
\item If $\Omega$ is exact, then the longitudinal KAM-cocycle of $\phi$ vanishes if and only
if $\Omega$ vanishes identically, i.e. $\phi$ is a geodesic flow;
\item If $\Omega$ is non-exact, then the longitudinal KAM-cocycle of $\phi$ vanishes
if and only if the Gaussian curvature is constant and $\Omega$ is a constant
multiple of the area form.
\end{enumerate}
}

\medskip

Item (1) was proved in \cite{P} using Aubry-Mather theory, but it was stated in a~different form.
The main result in \cite{P} asserts that if $\Omega$ is exact
and the Anosov splitting is of class $C^1$, then $\Omega$ must be zero (and this holds in any
dimension). The main result
of Foulon and Hasselblatt tells us that, for surfaces, the conditions of $C^1$ Anosov splitting
and longitudinal KAM-cocycle being a coboundary are equivalent.

The proof of item (2) in \cite{P2} for negative $K$ and small $\Omega$ was based on Fourier
analysis using the set up of V. Guillemin and D. Kazhdan in \cite{GK}.
Our approach here is based on establishing a Pestov identity for magnetic flows similar to
the ones in \cite{CS,DS} for geodesic flows.
Using this identity we will prove the following result which has independent interest:

\medskip

\noindent {\bf Theorem B.} {\it Let $M$ be a closed oriented surface and $\Omega$ 
an arbitrary smooth $2$-form. Suppose the magnetic flow $\phi$ of the pair
$(g,\Omega)$ is Anosov and let $X_{\la}$ be the vector field generating $\phi$.
If $G:M\to\mathbb R$ is any smooth
function and $\omega$ is any smooth $1$-form on $M$ such that
there is a smooth function $\varphi:SM\to \mathbb R$ for which
$G(x)+\omega_x(v)=X_{\la}(\varphi)$, then
$G$ is identically zero and $\omega$ is exact.

}

\medskip

Note that by the smooth Liv\v sic theorem \cite{LMM} saying that
$G(x)+\omega_x(v)=X_{\la}(\varphi)$ is equivalent to saying that
$G(x)+\omega_x(v)$ has zero integral over every closed magnetic geodesic.
Theorem A follows from Theorem B using the same methods as in \cite{P2}, so we will not repeat
the proof here. Instead we will consider a second application of Theorem B.

Given any closed 2-form $\Omega$, fix a constant $c\neq 0$ such that
the cohomology class of $c\,\Omega$ is an {\it integral class}, i.e.
$[c\,\Omega]\in H^2(M,\Z)=\Z$. 
Then there exists a principal circle bundle $\Pi:P\to M$ with
Euler class $[c\Omega]$. The bundle admits a connection 1-form $\alpha$
such that $d\alpha=-2\pi c\,\Pi^*\Omega$.
Recall that the holonomy function is a map $\log {\mbox{\rm hol}}_{\alpha}:Z_{1}(M)\to \re/\Z$,
where $Z_1(M)$ is the space of 1-cycles, such that for every 2-chain $f:\Sigma\to M$
we have
\[\log {\mbox{\rm hol}}_{\alpha}(\partial\Sigma)=-c\,\int_{\Sigma}f^*\Omega\;\;{\mbox{\rm mod}}\,1.\]

Let $\ga$ be a {\it closed} magnetic geodesic and let
$\ell(\ga)$ be its length. We define the {\it action} of $\ga$ as:
\[A(\ga):=\ell(\ga)-c^{-1}\,\log {\mbox{\rm hol}}_{\alpha}(\ga)\;{\mbox{\rm mod}}\,1.\]
We call the set ${\mathcal S}\subset \re/\Z$ of values $A(\ga)$ as $\ga$
ranges over all closed magnetic geodesics, the {\it action spectrum} of the magnetic flow.
If all the closed orbits of the magnetic flow $\phi$ are nondegenerate, then
${\mathcal S}$ is a countable set.

Suppose now that we vary
the connection 1-form $\alpha$. Let $\alpha_\tau$ be a smooth 1-parameter
family of connections for $\tau\in (-\varepsilon, \varepsilon)$ with
$\alpha_0=\alpha$ 
Then we can write $\alpha_{\tau}-\alpha=\Pi^*\beta_\tau$, where $\beta_\tau$
are smooth 1-forms on $M$.
The connection $\alpha_\tau$ has curvature form $-2\pi c\,\Omega+d\beta_\tau$.
If we let $\Omega_{\tau}=\Omega-\frac{1}{2\pi c}d\beta_\tau$
we get a magnetic flow $\phi^\tau$ and an action spectrum
${\mathcal S}_{\tau}$. If the magnetic flow $\phi$ is Anosov, then for
$\varepsilon$ small enough $\phi^\tau$ is Anosov for all 
$\tau\in (-\varepsilon, \varepsilon)$.

\medskip

\noindent {\bf Theorem C.} {\it Let $M$ be a closed oriented surface
endowed with a Riemannian metric $g$ and let $\Omega$ be a $2$-form.
Suppose the magnetic flow of the pair $(g,\Omega)$ is Anosov.
If ${\mathcal S}_{\tau}={\mathcal S}$ for all
$\tau $ sufficiently small, then the deformation is trivial, that is,
$\alpha_\tau=\alpha+\Pi^{*}dF_\tau$ and $\Omega_\tau=\Omega$, where $F_{\tau}$ are smooth
functions on $M$.}

\medskip

Theorem C and the results of V. Guillemin and A. Uribe in \cite{GU} give a version of infinitesimal
spectral rigidity for magnetic flows. In order to describe this rigidity we will assume that
$c=1$. This is really no restriction at all since the magnetic flows of $(g,\Omega)$ and
$(c^2\,g,c\,\Omega)$ are the same up to a constant time change.
For every positive integer $m$, let $L_m$ be the Hermitian line bundle with connection
over $M$ associated with $\Pi$ via the character $e^{i\theta}\mapsto e^{i m\theta}$ of $S^1$.
The metric on $M$, together with the connection on $L_m$ determine a Bochner-Laplace operator
acting on sections of $L_m$. For each $m$, let $\{\nu_{m,j}:\,j=1,2,\dots\}$ be the spectrum
of this operator. If we now vary the connection 1-form $\alpha$ as above we obtain 
eigenvalues $\nu_{m,j}^{\tau}$.

\medskip

\noindent {\bf Corollary.} {\it Let $M$ be a closed oriented surface
endowed with a Riemannian metric $g$ and let $\Omega$ be an integral $2$-form.
Suppose the magnetic flow of the pair $(g,\Omega)$ is Anosov.
If $\nu_{m,j}^\tau$ is independent of $\tau$ for all $m$ and $j$ $($i.e. the deformation
is isospectral$)$, then the deformation is trivial, that is,
$\alpha_\tau=\alpha+\Pi^{*}dF_\tau$ and $\Omega_\tau=\Omega$, where $F_{\tau}$ are smooth
functions on $M$.}

\medskip

Indeed, let us consider the periodic distribution
\[\Upsilon(s)=\sum_{m,j}\varphi\left(\sqrt{\nu_{m,j}+m^2}-m\sqrt{2}\right)e^{ims}\] 
where $\varphi$ is a Schwartz function on the real line. 
Theorem 6.9 in \cite{GU} asserts that the singularities of $\Upsilon$ are included in the
set of all $s\in\re$ for which $s/2\pi\,\mbox{\rm mod}\,1\in {\mathcal S}$. Moreover, each point
of the action spectrum arises as a singularity of $\Upsilon$ for an appropriate choice of $\varphi$.
The corollary is now an immediate consequence of Theorem~C.


\section{Preliminaries}

Let $M$ be a closed oriented surface, $SM$ the unit sphere bundle
and $\pi:SM\to M$ the canonical projection. The latter is in fact
a principal $S^{1}$-fibration and we let $V$ be the infinitesimal
generator of the action of $S^1$.

Given a unit vector $v\in T_{x}M$, we will denote by $iv$ the
unique unit vector orthogonal to $v$ such that $\{v,iv\}$ is an
oriented basis of $T_{x}M$. There are two basic 1-forms $\alpha$
and $\beta$ on $SM$ which are defined by the formulas:
\[\alpha_{(x,v)}(\xi):=\langle d_{(x,v)}\pi(\xi),v\rangle;\]
\[\beta_{(x,v)}(\xi):=\langle d_{(x,v)}\pi(\xi),iv\rangle.\]
The form $\alpha$ is the canonical contact form of $SM$ whose Reeb vector
field is the geodesic vector field $X$.
The volume form $\alpha\wedge d\alpha$ gives rise to the Liouville measure
$d\mu$ of $SM$.

A basic theorem in 2-dimensional Riemannian geometry asserts that
there exists a~unique 1-form $\psi$ on $SM$ (the connection form)
such that $\psi(V)=1$ and
\begin{align}
& d\alpha=\psi\wedge \beta\label{riem1}\\ & d\beta=-\psi\wedge
\alpha\label{riem2}\\ & d\psi=-(K\circ\pi)\,\alpha\wedge\beta
\label{riem3}
\end{align}
where $K$ is the Gaussian curvature of $M$. In fact, the form
$\psi$ is given by
\[\psi_{(x,v)}(\xi)=\left\langle \frac{DZ}{dt}(0),iv\right\rangle,\]
where $Z:(-\varepsilon,\epsilon)\to SM$ is any curve with
$Z(0)=(x,v)$ and $\dot{Z}(0)=\xi$ and $\frac{DZ}{dt}$ is the
covariant derivative of $Z$ along the curve $\pi\circ Z$.


For later use it is convenient to introduce the vector field $H$
uniquely defined by the conditions $\beta(H)=1$ and
$\alpha(H)=\psi(H)=0$. The vector fields $X,H$ and $V$ are dual to
$\alpha,\beta$ and $\psi$ and as a consequence of (\ref{riem1}--\ref{riem3}) they satisfy the commutation relations
\begin{equation}\label{comm}
[V,X]=H,\quad [V,H]=-X,\quad [X,H]=KV.
\end{equation}
Equations (\ref{riem1}--\ref{riem3}) also imply that the vector fields
$X,H$ and $V$ preserve the volume form $\alpha\wedge d\alpha$ and hence
the Liouville measure.

\section{Proof of Theorem B}

Henceforth $(M,g)$ is a closed oriented surface
and $X$, $H$, and $V$ are the same vector fields on $SM$
as in the previous section.

Let $\lambda$ be the smooth function on $M$ determined by $\Omega=\la\Omega_a$ and let
$$
X_\lambda=X+\lambda V
$$
be the generating vector field of the magnetic flow $\phi$ ($X_\la$ also preserves
Liouville measure).

From (\ref{comm}) we obtain:
$$
[V,X_\lambda]=H,\quad [V,H]=-X_\lambda+\lambda V,\quad 
[X_\lambda,H]=-\lambda X_\lambda+(K-H\lambda+\lambda^2)V.
$$
Note that
$$
H\lambda(x,v)=\langle \nabla \lambda(x), iv\rangle.
$$

\begin{Lemma}[Pestov's identity]\label{pestov}
For every smooth function $\varphi:SM\to\mathbb R$ we have
\begin{align*}
2H\varphi\cdot VX_\lambda\varphi
&=(X_\lambda\varphi)^2+(H\varphi)^2-(K-H\lambda+\lambda^2)(V\varphi)^2\\
&+X_\lambda(H\varphi\cdot V\varphi)
-H(X_\lambda\varphi\cdot V\varphi)
+V(X_\lambda\varphi\cdot H\varphi).
\end{align*}
\end{Lemma}

\begin{Remark}
{\rm A similar identity for the vector fields 
$X$, $H_\lambda:=H+\lambda V$ and $V$
was obtained in \cite[Lemma 2.1]{SU}}.
\end{Remark}

\begin{proof}[Proof of Lemma \ref{pestov}]
Using the commutation formulas, we deduce:
\begin{align*}
2H\varphi\cdot VX_\lambda&\varphi-V(H\varphi\cdot X_\lambda\varphi)\\
&=H\varphi\cdot VX_\lambda\varphi-VH\varphi\cdot X_\lambda\varphi\\
&=H\varphi\cdot(X_\lambda V\varphi+[V,X_\lambda]\varphi)
-X_\lambda\varphi\cdot(HV\varphi+[V,H]\varphi)\\
&=H\varphi\cdot(X_\lambda V\varphi+H\varphi)
-X_\lambda\varphi\cdot(HV\varphi-X_\lambda\varphi+\lambda V\varphi)\\
&=(X_\lambda\varphi)^2+(H\varphi)^2
+(X_\lambda V\varphi)(H\varphi)-(HV\varphi)(X_\lambda\varphi)
-\lambda X_\lambda\varphi\cdot V\varphi\\
&=(X_\lambda\varphi)^2+(H\varphi)^2
+X_\lambda(V\varphi\cdot H\varphi)
-H(V\varphi\cdot X_\lambda\varphi)
-[X_\lambda,H]\varphi\cdot V\varphi\\
&-\lambda X_\lambda\varphi\cdot V\varphi\\
&=(X_\lambda\varphi)^2+(H\varphi)^2
+X_\lambda(V\varphi\cdot H\varphi)
-H(V\varphi\cdot X_\lambda\varphi)\\
&-(K-H\lambda+\lambda^2)(V\varphi)^2
\end{align*}
which is equivalent to Pestov's identity.
\end{proof}

Integrating Pestov's identity over $SM$ against the Liouville measure $d\mu$,
we get
\begin{align}\label{id1}
2\int_{SM} H\varphi\cdot VX_\lambda\varphi\,d\mu
&=\int_{SM}(X_\lambda\varphi)^2\,d\mu+\int_{SM}(H\varphi)^2\,d\mu\\
&-\int_{SM}(K-H\lambda+\lambda^2)(V\varphi)^2\,d\mu.\notag
\end{align}

Let us derive one more integral identity.
Let $\varphi$ be again an arbitrary smooth function on $SM$.
By the commutation relations, we have
$$
X_\lambda V\varphi=VX_\lambda\varphi-H\varphi.
$$
Therefore,
$$
(X_\lambda V\varphi)^2=(VX_\lambda\varphi)^2+(H\varphi)^2
-2VX_\lambda\varphi\cdot H\varphi.
$$
Integrating, we obtain
\begin{equation}\label{id2}
\int_{SM}(X_\lambda V\varphi)^2\,d\mu
=\int_{SM}(VX_\lambda\varphi)^2\,d\mu
+\int_{SM}(H\varphi)^2\,d\mu
-2\int_{SM}VX_\lambda\varphi\cdot H\varphi\,d\mu.
\end{equation}

Subtracting (\ref{id1}) from (\ref{id2}), we
arrive at the final identity
\begin{align}\label{id}
\int_{SM}\left\{(X_\lambda V\varphi)^2-(K-H\lambda+\lambda^2)(V\varphi)^2\right\}&\,d\mu\\
&=\int_{SM}(VX_\lambda\varphi)^2\,d\mu-\int_{SM}(X_\lambda\varphi)^2\,d\mu.\notag
\end{align}

Let us now begin with the proof of Theorem B.
If $X_\lambda \varphi=G(x)+\omega_x(v)$, then it is easy to see that
the right-hand side of (\ref{id}) is nonpositive.
Indeed, since $\mu$ is invariant under $v\mapsto -v$ and $v\to iv$ we have
\[\int_{SM}\omega_{x}(v)\,d\mu=0\;\;\; \mbox{\rm and}\;\; \int_{SM}(\omega_{x}(v))^2\,d\mu=
\int_{SM}(\omega_{x}(iv))^2\,d\mu.\]
But $VX_\lambda\varphi=\omega_{x}(iv)$ and thus
\[\int_{SM}(VX_\lambda\varphi)^2\,d\mu-\int_{SM}(X_\lambda\varphi)^2\,d\mu=
\int_{SM}(G(x))^{2}\,d\mu\leq 0.\]

Setting $\psi=V\varphi$, we get
\begin{equation}\label{in}
\int_{SM}\left\{(X_\lambda \psi)^2
-(K-H\lambda+\lambda^2)\psi^2\right\}\,d\mu\le 0.
\end{equation}

We now show that this is possible if and only if $\psi=0$.
This would give $V\varphi=0$, which is equivalent
to the claim of the theorem.

\begin{Lemma}\label{ineq}
If $\phi$ is Anosov, then for every
closed magnetic geodesic $\gamma:[0,T]\to M$ and every
smooth function $z:[0,T]\to \mathbb R$ such that $z(0)=z(T)$
and $\dot z(0)=\dot z(T)$ we have
$$
I:=\int_0^T \left\{(\dot{z}^2
-(K-\langle\nabla\lambda,i\dot\gamma\rangle+\lambda^2)z^2\right\}\,dt\ge0
$$
with equality if and only if $z\equiv 0$.
\end{Lemma}

\begin{proof} 
Given $(x,v)\in SM$ and $\xi\in T_{(x,v)}TM$, let
\[J_{\xi}(t)=d_{(x,v)}(\pi\circ\phi_{t})(\xi).\]
We call $J_{\xi}$ a {\it magnetic Jacobi field} with initial
condition $\xi$. It was shown in \cite{PP1/2} that $J_{\xi}$
satisfies the following Jacobi equation:
\begin{equation}
\ddot{J_{\xi}} + R(\dot{\gamma},J_{\xi})\dot{\gamma} -[Y(\dot{J_{\xi}})+
(\nabla_{J_{\xi}}Y)(\dot{\gamma})] = 0,  \label{jacobi}
\end{equation}
where $\ga(t)=\pi\circ\phi_{t}(x,v)$, $R$ is the curvature tensor
of $g$ and $Y$ is determined by the equality $\Omega_{x}(u,v)=\langle Y_x(u),v\rangle$
for all $u,v\in T_{x}M$ and all $x\in M$.
Let us express $J_{\xi}$ as follows:
\[J_{\xi}(t)=x(t)\dot{\ga}(t)+y(t)i\dot{\ga}(t),\]
and suppose in addition that $\xi\in T_{(x,v)}SM$, which implies
\begin{equation}
g_{\ga}(\dot{J}_{\xi},\dot{\ga})=0.\label{tangente}
\end{equation}
A straightforward computation using (\ref{jacobi}) and (\ref{tangente}) shows that
$x$ and $y$ must satisfy the scalar equations:
\begin{align}
&\dot{x}=\la(\ga)\,y   \label{magJac1}\\
&\ddot{y}+\left[K(\ga)-\left\langle\nabla
\la(\ga),i\dot{\ga}\right\rangle+\la^{2}(\ga)\right]y=0.
\label{magJac2}
\end{align}
Let $E$ be the weak stable bundle of $\phi$.
Since for any $(x,v)\in SM$ the subspace $E$ does not intersect
the vertical subspace $\mbox{\rm Ker}\,d\pi_{(x,v)}$ \cite{PP,P3}, there exists a linear map $S(x,v):T_{x}M\to T_{x}M$
such that $E$ can be identified with the graph of $S$.
Let $u(x,v)$ be the trace of $S(x,v)$ and
let $J_{\eta}=x\dot{\ga}+yi\dot{\ga}$ be the Jacobi field with initial conditions $\eta=(iv,S(iv))\in E$.
Since $u(t)=\langle S(t)i\dot{\ga},i\dot{\ga}\rangle$ and $\dot{J}_{\eta}=SJ_{\eta}$
we see that
\begin{equation}
\dot{y}=u\,y
\label{chvar}
\end{equation}
Note that $y$ never vanishes. Given $z$ as in the hypothesis of the lemma, let $q$
be defined by the equation $z=qy$. Using equation (\ref{magJac2}) we have
\begin{align*}
I=-\int_{0}^{T}z(\ddot{z}+[K(\ga)-\left\langle\nabla
\la(\ga),i\dot{\ga}\right\rangle+\la^{2}(\ga)]z)\,dt&=
-\int_{0}^{T}q\,\frac{d}{dt}(\dot{q}y^2)\,dt\\
&=-[q\dot{q}y^2]_{0}^{T}+\int_{0}^{T}\dot{q}^2y^2\,dt.
\end{align*}
Using the periodicity properties of $z$ and (\ref{chvar}) we have
\[[q\dot{q}y^2]_{0}^{T}=[z\dot{q}y]_{0}^{T}=-[q\dot{y}z]_{0}^{T}=-[quyz]_{0}^{T}
=-[uz^2]_{0}^{T}\]
and the last term vanishes since $u$ is globally defined on $SM$.
We conclude that
\[I\geq 0\]
with equality if and only $\dot{q}\equiv 0$. Hence if $I=0$, $q$ must be a constant, which must be zero
since $y$ cannot be periodic in $T$.
\end{proof}

We continue now with the proof of Theorem B.
The last lemma, applied to the function $z=\psi(\gamma)$, yields
\begin{equation}\label{in_g}
\int_\gamma \left\{(X_\lambda \psi)^2
-(K-H\lambda+\lambda^2)\psi^2\right\}\,dt\ge 0
\end{equation}
for every closed magnetic geodesic $\gamma$.
Since the flow is Anosov, the invariant measures supported on 
closed orbits are dense in the space of all invariant measures
on $SM$.
Therefore, the above yields
$$
\int_{SM}\left\{(X_\lambda \psi)^2
-(K-H\lambda+\lambda^2)\psi^2\right\}\,d\mu\ge0.
$$
Combining this with (\ref{id}), we find that
\begin{equation}\label{in_sm}
\int_{SM}\left\{(X_\lambda \psi)^2
-(K-H\lambda+\lambda^2)\psi^2\right\}\,d\mu=0.
\end{equation}

By the non-negative version of the Liv\v sic theorem, proved independently by M. Pollicott and R. Sharp and by
A. Lopes and P. Thieullen (see \cite{LT, PS}),
we conclude from (\ref{in_g}) and (\ref{in_sm}) that
$$
\int_\gamma\left\{(X_\lambda \psi)^2
-(K-H\lambda+\lambda^2)\psi^2\right\}\,dt=0
$$
for every closed magnetic geodesic $\gamma$.
Applying again Lemma \ref{ineq}, we see that
$\psi$ vanishes on all closed magnetic geodesics.
Since the latter are dense in $SM$, the function $\psi$
vanishes on all of $SM$, as required.

\section{Proof of Theorem C}

We begin with a general easy lemma.
Given a smooth closed curve $\ga:[0,T]\to M$ and $k\in\re$ we define the {\it free time action of} $\ga$ as: 
\[A_{k}(\ga):=\frac{1}{2}\int_{0}^T|\dot{\ga}(t)|^{2}\,dt+kT-
c^{-1}\log {\mbox{\rm hol}}_{\alpha}(\ga)\;{\mbox{\rm mod}}\,1.\]
Recall that the {\it energy} is the function given by $E(x,v):=\frac{1}{2}|v|_{x}^{2}$.

\begin{Lemma} Let $\ga:[0,T]\to M$ be a closed magnetic geodesic with energy $k$.
Let $\ga_\tau: [0,T_{\tau}]\to M$, $\tau\in (-\varepsilon, \varepsilon)$,  be a smooth variation of $\ga$ by smooth closed curves.
Then
\[\frac{dA_{k}(\ga_{\tau})}{d\tau}(0)=0.\]
\label{criticalpoints}
\end{Lemma}

\begin{proof} The curves $\ga_\tau-\ga$ form a 1-cycle
which is the boundary of a 2-chain $\Sigma_\tau$. Then
\begin{equation}\label{pt}
c^{-1}\log {\mbox{\rm hol}}_{\alpha}(\ga_{\tau})-c^{-1}\log {\mbox{\rm hol}}_{\alpha}(\ga)=
-\int_{\Sigma_\tau}\Omega\;{\mbox{\rm mod}}\,1.
\end{equation}
If we let $W(t)$ be the variational vector field of $\ga_{\tau}$,
a straightforward calculation using that $\ga$ has energy $k$ and (\ref{pt}) 
shows that
\[\frac{dA_{k}}{d\tau}(0)=
-\int_{0}^{T}\left\langle \frac{D\dot{\ga}}{dt},W(t)\right\rangle\,dt
+\int_{0}^{T}\Omega(\dot{\ga}(t),W(t))\,dt.\]
Since $\ga$ is a magnetic geodesic, 
\[\frac{D\dot{\ga}}{dt}=Y_{\ga}(\dot{\ga})\]
and thus
\[\frac{dA_{k}}{d\tau}(0)=0.\]

\end{proof}

Let us assume now that we are under the hypotheses of Theorem C.

\begin{Lemma} Suppose ${\mathcal S}_{\tau}={\mathcal S}$ for all
$\tau\in (-\varepsilon, \varepsilon)$. Then
\[\int_{\ga_{\tau}}\frac{d\beta_\tau}{d\tau}=0\]
for every closed magnetic geodesic $\ga_\tau$ of $(g,\Omega_\tau)$.
\end{Lemma}

\begin{proof} In each nontrivial homotopy class we have a 1-parameter
family of closed magnetic geodesics $\ga_\tau$. 
Let 
\[a_{\tau}(\ga_{\tau}):=A^{\tau}_{1/2}(\ga_\tau)=
\ell(\ga)-c^{-1}\,\log {\mbox{\rm hol}}_{\alpha_{\tau}}(\ga)\;{\mbox{\rm mod}}\,1.\]
Since 
${\mathcal S}$ is countable and the map 
$(-\varepsilon, \varepsilon)\ni\tau\mapsto a_{\tau}(\ga_\tau)$ is continuous we have
\[a_{\tau}(\ga_\tau)=a_{\tau_{0}}(\ga_{\tau_{0}})\]
for all $\tau\in(-\varepsilon, \varepsilon)$. 
Since
\[\log {\mbox{\rm hol}}_{\alpha_{\tau}}(\sigma)=\log {\mbox{\rm hol}}_{\alpha_{\tau_{0}}}(\sigma)
+\frac{1}{2\pi}\int_{\sigma}(\beta_{\tau}-\beta_{\tau_{0}})\;{\mbox{\rm mod}}\,1\]
for all $\sigma$, we have:
\[a_{\tau}(\ga_\tau)=a_{\tau_{0}}(\ga_\tau)-\frac{1}{2\pi c}\int_{\ga_\tau}(\beta_\tau-\beta_{\tau_{0}})\;{\mbox{\rm mod}}\,1.\]
By Lemma \ref{criticalpoints}, the map $\tau\mapsto a_{\tau_{0}}(\ga_\tau)$
has a critical point at $\tau=\tau_{0}$, hence the last equality implies
\[\left.\frac{d}{d\tau}\right|_{\tau=\tau_{0}}\int_{\ga_{\tau}}(\beta_\tau-\beta_{\tau_{0}})=0\]
which is easily seen to imply
\[\int_{\ga_{\tau_{0}}}\left.\frac{d\beta_{\tau}}{d\tau}\right|_{\tau=\tau_{0}}=0.\]
\end{proof}

To complete the proof of Theorem C, observe that the previous lemma and 
Theorem~B imply that for each $\tau$, 
$\frac{d\beta_\tau}{d\tau}$ is exact. If we let $f_\tau$ be a primitive of $\frac{d\beta_\tau}{d\tau}$, then
\[F_\tau:=\int_{0}^\tau f_{s}\,ds\]
are the required functions.

\begin{Remark}{\rm The proofs of Theorem C and its corollary work
in any dimension provided that Theorem B holds in any dimension.
One only needs the cohomology class
$[\Omega]$ to be {\it rational}, i.e. there exists $\la\in\re$
such that $[\la\Omega]$ is an integral class.

Even if $[\Omega]$ is not rational, we can still attach to the magnetic flow an action spectrum
by considering a torus bundle $\T^r$ over $M$. The action spectrum is now
a subset of $\T^r$ and the same infinitesimal rigidity holds, provided that the
magnetic flow is Anosov.

The question of whether Theorem B extends to higher dimension is more delicate. We hope to discuss
these topics elsewhere, as well as the analogue of Theorem B for higher order tensors. 

}
\end{Remark}

\end{document}